\documentclass[11pt]{article}

\usepackage[margin=1in]{geometry}
\usepackage{amsfonts,amsmath,amsthm,amssymb}
\usepackage{graphicx}
\usepackage{xcolor}
\usepackage{booktabs}
\usepackage{multirow}
\usepackage{framed}
\usepackage{float}
\usepackage{algorithm}
\usepackage{algpseudocode}
\usepackage{url}
\usepackage[hidelinks]{hyperref}

\graphicspath{{./}{figures/}}

\providecommand{\bh}[1]{#1}

\newtheorem{thm}{Theorem}[section]

\newtheorem{lem}[thm]{Lemma}

\newtheorem{rem}[thm]{Remark}

\newenvironment{proofed}{\begin{proof}}{\end{proof}}

\title{fPINN-DeepONet: A Physics-Informed Operator Learning Framework for Multi-term Time-fractional Mixed Diffusion-wave Equations}

\author{
Binghang Lu\\
\small Department of Mathematics, Purdue University\\
\small West Lafayette, IN 47907, USA
\and
Zhao-peng Hao\\
\small Department of Mathematics, Purdue University\\
\small West Lafayette, IN 47907, USA
\and
Christian Moya\\
\small Department of Mathematics, Purdue University\\
\small West Lafayette, IN 47907, USA
\and
Guang Lin\thanks{Corresponding author. Email: \texttt{guanglin@purdue.edu}}\\
\small Department of Mathematics, Purdue University\\
\small School of Mechanical Engineering, Purdue University\\
\small West Lafayette, IN 47907, USA
}

\date{}

\begin{document}
\maketitle

\begin{abstract}
In this paper, we develop a physics-informed deep operator learning framework for solving multi-term time-fractional mixed diffusion-wave equations (TFMDWEs). We begin by deriving an $L_2$ approximation, which achieves first-order accuracy for the Caputo fractional derivative of order $\beta \in (1,2)$. Building upon this foundation, we propose the fPINN-DeepONet framework, a novel approach that integrates operator learning with the $L_2$ approximation to efficiently solve fractional partial differential equations (FPDEs). Our framework is successfully applied to both fixed and variable fractional-order PDEs, demonstrating the framework's versatility and broad applicability. To evaluate the performance of the proposed model, we conduct a series of numerical experiments that involve dynamically varying fractional orders in both space and time, as well as scenarios with noisy data. These results highlight the accuracy, robustness, and efficiency of the fPINN-DeepONet framework.
\end{abstract}

\noindent\textbf{Keywords:} $L_2$ approximation; operator learning; machine learning; data-driven scientific computing



\section{Introduction}

In this work, we focus on  operator learning methods for the multi-term time-fractional mixed diffusion-wave equations (TFMDWEs) in the following form~\cite{hao2016finite}:
\begin{eqnarray}\label{multi-dw}
&&\sum_{i=1}^sK_i\  _0^{C}D_t^{\alpha_i}u(\mathbf{x},t)= \Delta u(\mathbf{x},t) +f(\mathbf{x},t),\quad \mathbf{x}\in \Omega,~t\in (0, T],
\end{eqnarray}
where $\Omega$ is spatial domain, $\Delta$ is the Laplace operator,  $0< \alpha_1<\cdots \leq  \alpha_{i}\leq 1<\alpha_{i+1}< \cdots <\alpha_s <2 $ and $K_i>0,$ $i=1,\ldots,s.$ Here $ \  _0^{C}D_t^{\alpha_i}$ denotes the Caputo fractional derivative of order $\alpha_i,$ which will be defined later.

  As a natural extension of the single-term time-fractional partial differential equations (TFPDEs), e.g. sub-diffusion or diffusion-wave equations, the multi-term TFPDEs are expected to improve the modeling accuracy in depicting the anomalous diffusion process, successfully capturing power-law frequency dependence \cite{Kelly2008}, adequately modeling various types of viscoelastic damping \cite{Chen-Liu2012}. For example, in \cite{Schumer-Benson2003},  a two-term  mobile and immobile fractional-order diffusion model was proposed to model the total concentration in solute transport. The kinetic equation with two fractional derivatives of different orders appears quite natural when describing sub-diffusive motion in velocity fields \cite{Metzler1998}. The two-term time-fractional telegraph equations \cite{Orsingher2004} govern the iterated Brownian motion and the telegraph processes with Brownian
time.

    \bh{ {\bf Mathematical theories  on multi-term TFPDEs:}
 Based on an appropriate maximum principle, Luchko \cite{Luchko-2011} derived some priori
estimates for the solution and established its uniqueness using the Fourier
method and the method of separation of variables for the generalized multi-term time-fractional diffusion equation. Daftardar-Gejji and Bhalekar \cite{Daftardar-Gejji2008} considered a multi-term
time-fractional diffusion-wave equation (TFDWE) along with the homogeneous/nonhomogeneous boundary
conditions and have  solved this equation  using the method of separation of variables. Based
on the orthogonal polynomials of the Laguerre type, Stojanovic \cite{Stojanovic2011} found the solutions for the
diffusion-wave equation in one dimension with $n$-term time-fractional derivatives, whose orders
belong to the intervals (0, 1), (1, 2) and (0, 2), respectively. Using Luchko's Theorem and the
equivalent relationship between the Laplacian operator and the Riesz fractional derivative, Jiang
et al. \cite{Jiang2012} derived the analytical solutions for the multi-term time-space fractional advection-diffusion
equations. Subsequently, Ding et al. \cite{Ding2013} presented the analytical solutions
for the multi-term time-space fractional advection-diffusion equations with mixed boundary
conditions. Recently, Stojanovic has  analyzed regularity,  existence and uniqueness of the equation with nonlinear source  in  \cite{Stojanovic2013}.}\\
\bh{{\bf Numerical approximations of multi-term TFPDEs:}  For the multi-term TFPDEs  whose orders, $\alpha_i$, belong to $(0,1),$ very recently, Jin et al. \cite{Jin2015} developed a fully discrete scheme based  on the standard Galerkin
finite element method in space  combining with a finite difference discretization of the time-fractional derivatives.   For the orders, $\alpha_i$, lying in  $(1,2),$ Chen. et al.  \cite{Chen-Liu2012} presented a finite difference scheme and gave its analysis, following the idea of the method of order reduction  proposed by Sun and Wu \cite{Sun2006}. In addition, an extension of multi-term TFPDEs, distributed order TFPDEs have also been considered (see e.g.  \cite{Ye2014,Morgado2014}). However, traditional methods (e.g., finite difference, finite element, spectral methods) require significant computational resources due to the non-local nature of fractional derivatives and high-dimensional FPDEs (e.g., multi-term or variable-order FPDEs) suffer from exponential growth in computational cost, making grid-based methods inefficient.} \\
\bh{{\bf Deep learning methods for solving fractional partial differential equations:} Xingjian Xu et al.\cite{xu2022discovery} employed a deep neural network (DNN) combined with the $L_1$ approximation to investigate subdiffusion problems. Guofei Pang et al.\cite{pang2019fpinns} introduced fractional physics-informed neural networks (fPINNs) to solve space-time fractional advection-diffusion equations (fractional ADEs). Fang et al.\cite{fang2023explore} proposed the use of DNNs for solving FPDEs and their associated inverse problems. Guo\cite{guo2022monte} developed Monte Carlo-based fPINNs to address both forward and inverse FPDEs. However, the integration of the $L_2$ approximation with deep operator learning methods for solving FPDEs involving multi-term and variable fractional orders remains an unexplored area. }

In this paper, a fPINN-DeepONet framework, a physics-informed operator learning method, is proposed to solve multi-term TFPDEs. Specifically, a fPINN-DeepONet, $G_\theta$, is created to approximate the solution operator~$G^\dagger: u \mapsto s$. Here $u$ can be the variable fractional order $\alpha_i$ to ensure the model to generalize over different fractional orders or the input function $[u(x,0),f(x,t)]$. Our numerical results demonstrate that the proposed framework achieves accurate prediction in various scenarios, including forward and inverse problems, high-dimensional fractional orders, discontinuous orders, and in the presence of noisy data.

In addition to the proposed fPINN-DeepONet framework, several recent studies have explored related directions in scientific machine learning, physics-informed learning, optimization, and modern machine learning applications. For physics-informed scientific computing, iterative PINN frameworks combined with ensemble Kalman filtering have been developed to improve solution accuracy and training robustness~\cite{lu2025ipinner}. Evolutionary multi-objective optimization has also been incorporated into replica-exchange-based physics-informed neural operator learning to enhance the efficiency and reliability of neural operator training~\cite{lu2025morephynet}. More recent efforts have investigated adaptive momentum feedback-linearization methods for hard-constrained PINN optimization~\cite{lu2026adamflip}, spectral-guided Muon optimization for scientific machine learning~\cite{lu2026spectralguidance}, and neural-operator-based frameworks for infinite-dimensional functional nonlinear proper orthogonal decomposition~\cite{mou2026neuralpod}. Beyond scientific machine learning, related work has also examined  object detection for mask-wearing recognition~\cite{liu2020yolov5mask}, spectral orthogonal gradient projection for continual learning in large language models~\cite{lu2026muonogd}, and safety mitigation for large reasoning models through adaptive multi-principle steering~\cite{li2026chainrisk}. These studies provide broader context for the development of reliable, efficient, and generalizable learning-based methods for complex scientific and engineering problems.

The remainder of this paper is organized as follows. In Section 2, we present a first-order $L_2$ approximation for the Caputo fractional derivative and derive the corresponding compact difference scheme. Section 3 provides an introduction to DeepONet and operator learning. We then describe the proposed fPINN-DeepONet framework and explain the construction of the loss function, including the implementation of the $L_2$ approximation for computing the physics-based loss term. A pseudocode outlining the training procedure is provided at the end of the section. In Section 4, we present experimental results on both fixed and variable fractional-order partial differential equations to evaluate the performance of fPINN-DeepONet. We also assess the robustness of our method in the presence of noisy data. Section 5 discusses hyperparameter tuning to optimize model performance. Finally, we conclude the paper in Section 6.

\section{Numerical method for the one-dimensional two-term time-fractional  mixed diffusion-wave equation}
\label{2.numerical method}
\setcounter{equation}{0}

For the sake of simplicity but without loss of  generality,  we  consider the one-dimensional two-term TFMDWEs as follows:
\begin{eqnarray}
&&K_1\  _0^{C}D_t^{\alpha}u(x,t)+K_2\  _0^{C}D_t^{\beta}u(x,t)= \partial_x^{2}u(x,t) +f(x,t),\quad x\in \Omega,~0<t\leq T, \label{problem-1}\\
&&u(x,0)=\phi_0(x),\quad u_t(x,0)=\phi_1(x),\quad x\in \bar{\Omega},\label{problem-2} \\
 &&u(0,t)=\varphi_0(t),\quad u(L,t)=\varphi_1(t),\quad 0< t\leq T ,\label{problem-3}
\end{eqnarray}
where $\Omega=(0,L),$ $K_1,K_2> 0,$   and  $0< \alpha< 1 <\beta< 2.$

In the following analysis of the proposed numerical method, we assume that the problem \eqref{problem-1}-\eqref{problem-3} has a unique and sufficiently smooth solution.
Here, it should be noted that Alikhanov showed the stability of the solution by an a priori estimate for the subdiffusion and diffusion wave equations in \cite{Alikhanov2010}.  Following his idea, the stability for the solution of problem \eqref{problem-1}-\eqref{problem-3} can be investigated by the energy method.
\subsection{Notations and lemmas}

Take an integer $N.$ Let $\Omega_\tau\equiv \{ t_n~|~0\leq n\leq N\}$ be a uniform mesh of the interval $[0,T],$ where $t_n=n\tau,~0\leq n\leq N$ with $\tau=T/N.$   Two lemmas below are needed. The first is a well-known $L1$ approximation.
 \begin{lem}(See \cite{Sun2006})\label{dis-lem-sun-1}
For $\alpha \in (0,1), $ suppose $v(t) \in C^2[0,t_n]$ and denote $v^n=v(t_n),$
\begin{equation} 
a_0^{(\alpha)}=\frac{1}{\Gamma{(2-\alpha)}}, \quad  a_n^{(\alpha)}=
\frac{1}{\Gamma{(2-\alpha)}}[(n+1)^{1-\alpha}-n^{1-\alpha}],\qquad n\geq 1, \label{coefficients}\end{equation}
 and define the backward difference quotient operator
\begin{equation}\label{eq:def-al}
\delta_t^{\alpha}v^n= \frac{1}{ \tau^{\alpha} }\biggl[a_0^{(\alpha)}v^n- \sum_{k=1}^{n-1} (a_{n-k-1}^{(\alpha)}-a_{n-k}^{(\alpha)} )v^k -a_{n-1}^{(\alpha)}v^0\biggl],\quad n \geq 1.
\end{equation}
It holds that
\begin{equation}
_0^{C}D_t^{\alpha}v(t_n)=\delta_t^{\alpha}v^n +R_1[v(t_n)],
\end{equation}
where
$$|R_1(v(t_n))|\leq c_{\alpha}\max_{0\leq t\leq t_n}|v^{\prime\prime}(t)|\tau^{2-\alpha}. $$
\end{lem}

 In the above lemma $c_{\alpha}$ is a bounded constant, which depends on $\alpha$, but is independent of step size \bh{$\tau$.} 
For the approximation of the derivative of fractional order $\beta \in(1,2),$ we have the following result.
\begin{lem}\label{dis-lemma-2}
For $\beta \in(1,2), $  suppose  $v(t)\in C^{2}[0,t_n]\cap C^{3}(0,t_n], $ and $|v^{\prime\prime\prime}(t)|\in L_1(0,t_n].$  Denote $v^n=v(t_n),$
\begin{equation}
b_0^{(\beta)}:=a_{0}^{(\beta-1)}=\frac{1}{\Gamma{(3-\beta)}},\quad b_n^{(\beta)}:=a_n^{(\beta-1)}=
\frac{1}{\Gamma{(3-\beta)}}[(n+1)^{2-\beta}-n^{2-\beta}],\quad n\geq 1 \label{coefficients-2}\end{equation}
and
\begin{eqnarray}\label{eq:def-be}
&&\Delta_t^{\beta}v^n=
\frac{1}{\tau^{\beta}}\biggl[  \sum_{k=2}^{n} b_{n-k}^{(\beta)}(v^k-2v^{k-1}+v^{k-2})+ 2b_{n-1}^{(\beta)}(v^1-v^0) \biggl].
\end{eqnarray}  It holds that
\begin{eqnarray} \label{denotion-be}
&& _0^CD_t^{\beta}v(t_n)=\Delta_t^{\beta}v^n- 2  \frac{b_{n-1}^{(\beta)}}{\tau^{\beta-1}}  v^{\prime}(t_0)  +R_2[v(t_n)],\quad n\geq 1,
\end{eqnarray}
where
$$ |R_2[v(t_n)]|\leq \frac{9t_n^{2-\beta}}{\Gamma{(3-\beta)}} \max_{0\leq k\leq n-1}\int_{0}^{1}|v^{\prime\prime\prime}(t_{k}+\theta \tau)|d\theta \cdot \tau.$$
\end{lem}
\begin{proofed}
Note that
\begin{eqnarray}
&&_0^CD_t^{\beta}v(t_n)=\frac{1}{\Gamma{(2-\beta)} } \sum_{k=2}^n \int_{t_{k-1}}^{t_k} \frac{v^{\prime\prime}(s)}{(t_n-s)^{\beta-1}}ds +\frac{1}{\Gamma{(2-\beta)} } \int_{0}^{t_1} \frac{v^{\prime\prime}(s)}{(t_n-s)^{\beta-1}}ds\nonumber \\
&&   \qquad \qquad ~ = \frac{1}{\Gamma{(2-\beta)} } \sum_{k=2}^n \int_{t_{k-1}}^{t_k} \frac{\delta_t^2 v^{k-1}}{(t_n-s)^{\beta-1}}ds +\frac{1}{\Gamma{(2-\beta)} } \int_{t_{0}}^{t_1} \frac{\Delta_t^2v^{0}}{(t_n-s)^{\beta-1}}ds+R_2[v(t_n)] \nonumber \\
&&   \qquad \qquad~=\Delta_t^{\beta}v^n- 2  \frac{b_{n-1}^{(\beta)}}{\tau^{\beta-1}}  v^{\prime}(t_0)  +R_2[v(t_n)], \label{split-1}
\end{eqnarray}
where
$$\delta_t^2v^{k-1}=\frac{v(t_k)-2v(t_{k-1})+v(t_{k-2}) }{\tau^2},\quad \Delta_t^2v^{0}=2\frac{v(t_1)-v(t_{0})-\tau v^{\prime}(t_{0}) }{\tau^2},$$
and
\begin{eqnarray}
&&R_2[v(t_n)]= \frac{1}{\Gamma{(2-\beta)} } \sum_{k=2}^n \int_{t_{k-1}}^{t_k} \frac{v^{\prime\prime}(s)-\delta_t^2 v^{k-1}}{(t_n-s)^{\beta-1}}ds  +\frac{1}{\Gamma{(2-\beta)} }  \int_{t_{0}}^{t_1} \frac{v^{\prime\prime}(s)-\Delta_t^2 v^{0}}{(t_n-s)^{\beta-1}}ds. \label{truncation}
\end{eqnarray}
By  Taylor expansion with the integral remaining term,  we find for $s\in(t_{k-1},t_{k})$ that
\begin{eqnarray}
&& |v^{\prime\prime}(s)-\delta_t^2 v^{k-1}| = \biggl|\int^s_{t_{k-2}}\biggl(1-\frac{1}{2\tau^2}[(t_k-t)^2+(t_k-s)^2]\biggl) v^{\prime\prime\prime}(t)dt-\int_s^{t_{k}}\frac{1}{2\tau^2}(t_{k}-t)^2 v^{\prime\prime\prime}(t)dt\nonumber\\
&&\qquad\qquad\qquad\qquad +\int_{t_{k-2}}^{t_{k-1}}\frac{1}{\tau^2}(t_{k-1}-t)^2 v^{\prime\prime\prime}(t)dt \biggl|\nonumber\\
&& \qquad \qquad\qquad\quad \leq \frac{7}{2}\int^s_{t_{k-2}} |v^{\prime\prime\prime}(t)|dt+\frac{1}{2}\int_s^{t_{k}} |v^{\prime\prime\prime}(t)|dt+\int_{t_{k-2}}^{t_{k-1}} |v^{\prime\prime\prime}(t)|dt\nonumber\\
&& \qquad \qquad\qquad\quad \leq 9\max_{1\leq k\leq n}\int_{t_{k-1}}^{t_{k}} |v^{\prime\prime\prime}(t)|dt =9 \max_{0\leq k\leq n-1}\int_{0}^{1}|v^{\prime\prime\prime}(t_{k}+\theta \tau)|d\theta\cdot \tau . \label{truncation-1}
\end{eqnarray}
Similarly we can obtain
\begin{eqnarray}
 |v^{\prime\prime}(s)-\Delta_t^2 v^{0}|\leq \max_{0\leq k\leq n-1}\int_{0}^{1}|v^{\prime\prime\prime}(t_{k}+\theta \tau)|d\theta\cdot \tau. \label{truncation-2}
\end{eqnarray}
Combining  \eqref{truncation-1}-\eqref{truncation-2} with \eqref{truncation}, we get
\begin{eqnarray}
&& |R_2[v(t_n)]| \leq  \max_{0\leq k\leq n-1}\int_{0}^{1}|v^{\prime\prime\prime}(t_{k}+\theta \tau)|d\theta \cdot \frac{9}{\Gamma{(2-\beta)} } \sum_{k=1}^n \int_{t_{k-1}}^{t_k} \frac{1}{(t_n-s)^{\beta-1}}ds \cdot \tau
\nonumber \\
&&\qquad\qquad~
 = \frac{9t_n^{2-\beta}}{\Gamma{(3-\beta)}} \max_{0\leq k\leq n-1}\int_{0}^{1}|v^{\prime\prime\prime}(t_{k}+\theta \tau)|d\theta \cdot \tau.
\end{eqnarray}
This completes the proof.
\end{proofed}

\bh{\subsection{Derivation of the compact finite difference scheme}
We are now in a position to derive the compact finite difference scheme.
Define grid functions below  
$$U_i^n=u(x_i,t_n), \quad (U_t)_i^0= u_t(x_i,t_0)\quad F_i^n=f(x_i,t_n),\quad 0\leq i\leq M,\quad 0\leq n\leq N. $$
Considering Equation \eqref{problem-1} at grid points  $ (x_i,t_n),$, we have
\begin{eqnarray}\label{eq-5}
&&K_1\  _0^{C}D_t^{\alpha}u(x_i,t_n)+K_2\  _0^{C}D_t^{\beta}u(x_i,t_n)=  \partial_x^2 u(x_i,t_n) +f(x_i,t_n),\quad 0\leq i\leq M,\quad 0\leq n\leq N.
\end{eqnarray}
By Lemmas \ref{dis-lem-sun-1} and \ref{dis-lemma-2}, we have
\begin{eqnarray}\label{eq-6}
&&K_1 \delta_t^{\alpha}U_i^n  +K_2 \Delta_t^{\beta} U_i^n= \partial_x^2u (x_i,t_n)+F_i^n+\frac{2K_2b_{n-1}^{(\beta)}}{\tau^{\beta-1}} u_t(x_i,0)+K_1R_t^1[u(x_i,t_n)]+K_2R_t^2[u(x_i,t_n)], \nonumber\\
&&\quad 0\leq i\leq M,\quad 1\leq n\leq N.
\end{eqnarray}
Here $R_1$ and $R_2$ are similarly  defined as that in Lemmas \ref{dis-lem-sun-1} and \ref{dis-lemma-2} respectively.
For the spatial  discretization, acting the average operator $\mathcal{A}_x$ that
\begin{eqnarray}\label{eq-7}
&&K_1 \mathcal{A}_x\delta_t^{\alpha}U_i^n  +K_2 \mathcal{A}_x\Delta_t^{\beta} U_i^n=  \delta_x^2 U_i^n+\mathcal{A}_xF_i^n+   \frac{2K_2b_{n-1}^{(\beta)}}{\tau^{\beta-1}} \mathcal{A}_x (U_t)_i^0+R_i^n, \quad 1\leq i\leq M-1,\quad 1\leq n\leq N.\nonumber \\
\end{eqnarray}
where
$$R_i^n=K_1\mathcal{A}_xR_t^1[u(x_i,t_n)]+K_2\mathcal{A}_xR_t^2[u(x_i,t_n)]+R_x^3[u(x_i,t_n)].$$
 and $$R_x^3[u(x_i,t_n)]=\frac{h^4}{360}\int_0^1[\partial_x^6u(x_i-sh,t_n)+\partial_x^6u(x_i+sh,t_n)]\theta (s)ds.$$
 It follows from Lemmas  \ref{dis-lem-sun-1} and \ref{dis-lemma-2}  that  there exists a constant  $C_u,$ which depends on the regularity of the solution $u(x,t)$ and the parameters $\alpha$ and $\beta$, but is independent of the step size $h$ and $\tau,$
such that
\begin{eqnarray}
&&|R_i^n|\leq C_u(\tau+h^4),\quad 1\leq i\leq M-1,~ 1\leq n\leq N.\label{truncation-error-1}
\end{eqnarray}
Noticing the initial-boundary conditions,
\begin{eqnarray}
&&U_i^0=\phi_0(x_i) , \quad  (U_t)_i^0=\phi_1(x_i),\quad  0\leq i\leq M , \label{eq-8}   \\
&& U_0^n=\varphi_0(x_0),\quad U_M^n=\varphi_1(x_M),\quad 1\leq n\leq N ,     \label{eq-9}
\end{eqnarray}
 omitting the small terms $R_i^n $ and
denoting by $u_i^n$ the numerical approximation of $U_i^n,$  we get the compact finite difference scheme,
\begin{eqnarray}
&&K_1 \mathcal{A}_x\delta_t^{\alpha}u_i^n  +K_2 \mathcal{A}_x\Delta_t^{\beta} u_i^n=  \delta_x^2 u_i^n+\mathcal{A}_xF_i^n+   \frac{2K_2b_{n-1}^{(\beta)}}{\tau^{\beta-1}} \mathcal{A}_x \phi_i^1,
\quad
1\leq i\leq M-1,~1\leq n\leq  N,\quad \label{b9}\\
&&u_i^0=\phi_0(x_i) , \quad  0\leq i\leq M , \label{b10}   \\
&&  u_0^n=\varphi_0(t_n),\quad u_M^n=\varphi_1(t_n),\quad 1\leq n\leq N.  \label{b11}
\end{eqnarray}
\begin{rem}
When $\alpha=1,$ $\beta=2,$ we get the following three time-level backward difference scheme
\begin{eqnarray}
&&K_1 \mathcal{A}_x \biggl(\frac{u_i^n-u_i^{n-1}}{\tau}\biggl)  +K_2 \mathcal{A}_x \biggl( \frac{u_i^n-2u_i^{n-1}+u_i^{n-2}}{\tau^2} \biggl)=  \delta_x^2 u_i^n+\mathcal{A}_xF_i^n,
 \nonumber\\
&&
\quad
1\leq i\leq M-1,~2\leq n\leq  N,\quad \label{b99}\\
&&K_1 \mathcal{A}_x \biggl(\frac{u_i^1-u_i^{0}}{\tau} \biggl)  + 2K_2  \mathcal{A}_x \biggl(\frac{u_i^1-u_i^{0}-\tau \phi_i}{\tau^2} \biggl) =  \delta_x^2 u_i^n+\mathcal{A}_xF_i^n,
\quad
1\leq i\leq M-1,\quad \label{b112}\\
&&u_i^0=\phi_0(x_i) , \quad  0\leq i\leq M , \label{b100}   \\
&&  u_0^n=\varphi_0(t_n),\quad u_M^n=\varphi_1(t_n),\quad 1\leq n\leq N.  \label{b111}
\end{eqnarray}
\end{rem}
}

\section{fPINN-DeepONet framework}
Operator learning involves training neural networks to approximate operators, which are mappings between infinite-dimensional function spaces. Instead of learning a function that maps inputs to outputs, DeepONet learns an operator, enabling it to handle functional dependencies effectively.
\begin{figure}[H]
     \centering
       \includegraphics[width=1.0\textwidth]{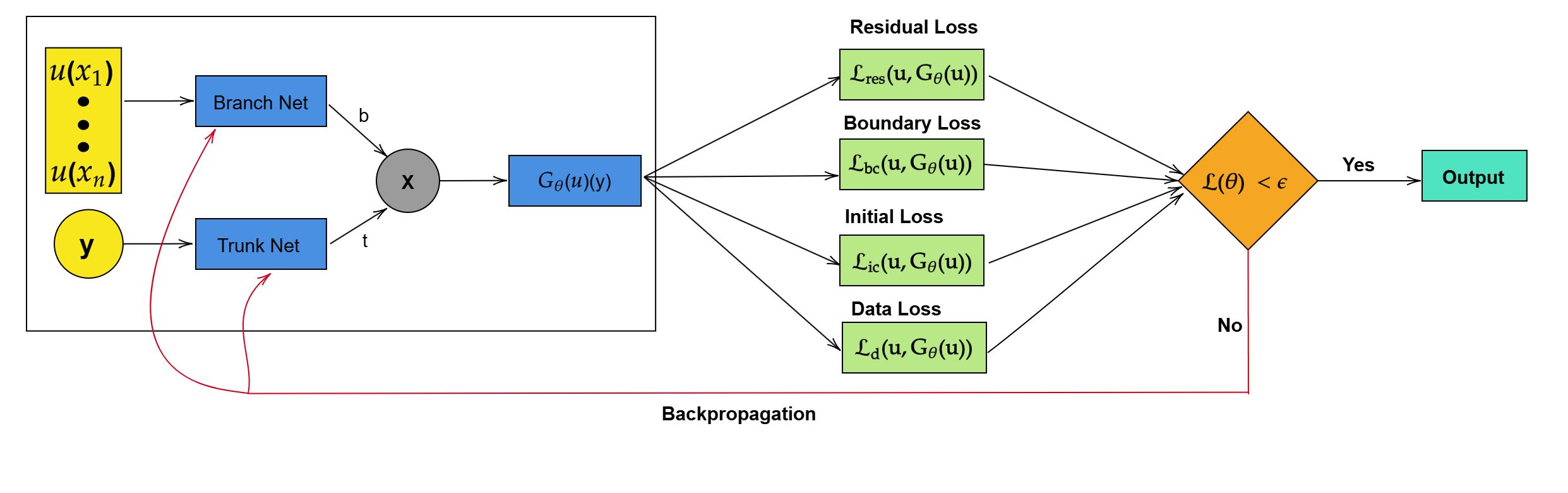}
     \hfill
     \caption{fPINN-DeepONet structure: The Branch Net takes the function $u$ as input, while the Trunk Net receives $y$ from the domain of $G(u)$. The output is given by $G_{\theta}(u)(y)$, and the total loss comprises the residual loss, boundary loss, initial loss, and data loss.}
    \label{fig:DeepONet framework}
\end{figure}
The Branch network processes the input function $u$ by discretizing it at a set of predefined sensor points. Let $(x_1, \dots, x_m)$ denote the collection of $m$ interpolation points, commonly referred to as sensors in \cite{lu2019deeponet}, which are used to encode the input function $u$. Specifically, we define the discretized input as $u_m := (u(x_1), \dots, u(x_m))$. The Branch network maps this discretized input $u_m$ to a vector of trainable coefficients $b \in \mathbb{R}^q$.

Similarly, the Trunk network processes the query location $y \in Y$, mapping it to a vector representing a collection of trainable basis functions \cite{moya2022fed}. The final output is a real number, obtained by combining the outputs of the Branch and Trunk networks via a dot product: the Branch output $b$ and the Trunk output $t$ are merged as following:

\begin{eqnarray}
G_{\theta}(u)(y) = \sum_{i=1}^q b_it_i
\end{eqnarray}

In this work, we propose a fractional-order Physics-informed DeepONet (fPINN-DeepONet)~$G_\theta$ framework, \bh{constructed} on the concept of operator learning, to approximate the solution operator~$G^\dagger: \alpha(t) \mapsto s$. Here, $\alpha(t)$ denotes the time-dependent fractional order in the FPDE, and $s$ is the corresponding solution of interest. As illustrated in Figure~\ref{fig:DeepONet framework}, the input $y$ to the trunk net lies in the domain of $G(u)$ and consists of spatial-temporal coordinates $x$ and $t$. The input to the branch net is the function $\alpha(x,t)$. 

For fixed fractional-order FPDEs, the branch net input can instead be the function $u(x,t)$, the right-hand side $f$, or the initial condition. The fPINN-DeepONet framework incorporates the governing physical laws of FPDEs~\cite{lu2023nsga} into the loss function, which is composed of two parts: the physical loss $\mathcal{L}_{phy}$ and the data loss $\mathcal{L}_{data}$. The physical loss includes the residual loss, initial condition loss, and boundary condition loss:

\begin{eqnarray} 
\mathcal{L}_{phy} =  \mathcal{L}_{res} + \mathcal{L}_{ic} + \mathcal{L}_{bc}
\end{eqnarray}
The residual loss can be calculated as 
\begin{eqnarray}
  \mathcal{L}_{res} = \frac{1}{N} \sum_{i=1}^{N} \left( K_1\  _0^{C}D_t^{\alpha}u(x_i,t_i) + K_2\  _0^{C}D_t^{\beta}u(x_i,t_i) - \partial_x^{2}u(x_i,t_i) - f(x_i,t_i) \right)^2 
\end{eqnarray}
\bh{As the partial derivative term $\partial_x^{2}u(x_i,t_i)$ can be solved by automatic differentiation with PyTorch~\cite{paszke2017automatic}. The source term $f(x,t)$ is known by solving the FPDE. The key part is to solve the first, and the second terms $K_1\  _0^{C}D_t^{\alpha}u(x_i,t_i) + K_2\  _0^{C}D_t^{\beta}u(x_i,t_i)$ using the predicted value from the neural network. To address this, we apply the $L_2$ approximation in Eqns.~\ref{coefficients} and~\ref{coefficients-2}, to solve these two terms. }\\
Data loss is calculated by measuring the mean squared error (MSE) between the predicted solution $u_\text{approx}$ and the true solution $u_\text{true}$.
\begin{eqnarray}
\mathcal{L}_{\text{data}} = \frac{1}{N} \sum_{i=1}^{N} \left| u_\text{approx}(x_i,t_i) - u_\text{true}(x_i,t_i) \right|^2
\end{eqnarray}
The total loss function is defined by:
\begin{eqnarray} \label{fPINN-DeepONet loss function}
\mathcal{L} =   \mathcal{L}_{phy} + \mathcal{L}_{data} 
\end{eqnarray}
The relative error $\| \cdot \|_r$ is measured by:
\begin{eqnarray}
\| \cdot \|_r = \frac{\sqrt{\sum_{i=1}^{N} \left( u_{\text{approx}}(x_i,t_i) - u_{\text{true}}(x_i,t_i) \right)^2}}{\sqrt{\sum_{i=1}^{N} \left( u_{\text{true}}(x_i,t_i) \right)^2}} 
\label{mse_error}
\end{eqnarray}
Here we employ the proposed fPINN-DeepONet framework to solve the FPDEs with fixed and variable fractional orders.
\newpage

\begin{algorithm}[H]
\caption{Solving fixed fractional-order FPDE using fPINN-DeepONet}
\label{alg:fixed-fpde}
\begin{algorithmic}[1]
\Require Epoch $I$; fixed fractional orders $\alpha$, $\beta$
\Ensure $u_{\mathrm{approx}}$
\State Generate mesh data $t$ with step size $\tau_t=100$ on the interval $[0,1]$ and $x$ with step size $\tau_x=100$ on the interval $[0,\pi]$.
\State Calculate the real solution $u_{\mathrm{true}}$ with $t$, $x$, $\alpha$, and $\beta$.
\State Calculate the right-hand side $f_{\mathrm{true}}(x,t)$ with $u_{\mathrm{true}}$, $x$, and $t$.
\State Construct the input function $f(x,t)$ for the branch net and $[x,t]$ for the trunk net.
\State Initialize the weights and biases in fPINN-DeepONet.
\For{$i=1$ to $I$}
    \State Construct fPINN-DeepONet for the training data and predict the output $u_{\mathrm{approx}}$.
    \State Use the $L_2$ approximation to solve $u_{\mathrm{approx}}$ for $f_{\mathrm{approx}}$ using Eqs.~\eqref{dis-lem-sun-1} and \eqref{dis-lemma-2}.
    \State Update the fPINN-DeepONet by minimizing the loss function in Eq.~\eqref{fPINN-DeepONet loss function}.
\EndFor
\State Calculate the $L_2$ testing error using Eq.~\eqref{mse_error}.
\end{algorithmic}
\end{algorithm}

The algorithm to solve the fixed fractional order FPDE is explained in detail in Algorithm 1. The input data $x,t$ are generated by grid sampling in the domain $x\in [0,\pi]$ and $t \in [0,1]$ for the trunk net, and the source term $f(x,t)$ is the input for the branch net. Since our proposed model is physics-informed, the initial condition $u(x,0)$ is considered in the loss term, there is no need to incorporate it in the branch net. The fPINN-DeepONet was used to approximate the solution, $u_{approx}$. Using the $L_2$ approximation method we proposed in \bh{Section~\ref{dis-lem-sun-1}}, we can calculate the right-hand side $f_{approx}$. The model was updated by minimizing the loss function ~\ref{fPINN-DeepONet loss function}.

\begin{algorithm}[H]
\caption{Solving variable fractional-order FPDE using fPINN-DeepONet}
\label{alg:variable-fpde}
\begin{algorithmic}[1]
\Require Epoch $I$
\Ensure $u_{\mathrm{approx}}$
\State Generate mesh data $t$ with step size $\tau_t=100$ on the interval $[0,1]$ and $x$ with step size $\tau_x=100$ on the interval $[0,\pi]$.
\State Generate variable fractional orders $\alpha(x,t)$ and $\beta(x,t)$.
\State Calculate the real solution $u_{\mathrm{true}}$ with $t$, $x$, $\alpha(x,t)$, and $\beta(x,t)$.
\State Calculate the right-hand side $f_{\mathrm{true}}(x,t)$ with $u_{\mathrm{true}}$, $x$, and $t$.
\State Construct the input functions $\alpha(x,t)$ and $\beta(x,t)$ for the branch net and $[x,t]$ for the trunk net.
\State Initialize the weights and biases in fPINN-DeepONet.
\For{$i=1$ to $I$}
    \State Construct fPINN-DeepONet for the training data and predict the output $u_{\mathrm{approx}}$.
    \State Use the $L_2$ approximation to solve $u_{\mathrm{approx}}$ for $f_{\mathrm{approx}}$ using Eqs.~\eqref{dis-lem-sun-1} and \eqref{dis-lemma-2}.
    \State Update the fPINN-DeepONet by minimizing the loss function in Eq.~\eqref{fPINN-DeepONet loss function}.
\EndFor
\State Calculate the $L_2$ testing error using Eq.~\eqref{mse_error}.
\end{algorithmic}
\end{algorithm}

In Algorithm 2, we extend the fractional orders, $\alpha$ and $\beta$ to be functions of both $x$ and $t$ denoted as $\alpha(x,t)$, $\beta(x,t)$, and employ the $L_2$ approximation method accordingly to address this extension.  Here, the fractional-order functions, $\alpha(x,t)$, $\beta(x,t)$ were used as the input functions for the branch net and the spatiotemporal data $[x,t]$ were used as the input data for the trunk net. The model is updated by minimizing the loss function~\ref{fPINN-DeepONet loss function}.

\section{Experimental Results}
\label{7.experimental results}
In this section, we present numerical experiments using the proposed fPINN-DeepONet to solve FODEs and FPDEs with fixed fractional order $\alpha_i$ within (0,2) are presented. We also illustrate numerical results for solving FPDEs with variable fractional order  $\alpha(t)$ and two-dimensional fractional order $\alpha(x,t)$. Additionally, we examine a case involving high spatial dimensions. The robustness of our model is tested by introducing different levels of Gaussian noise. Finally, we demonstrate the prediction results for solving FPDEs with a discontinuous fractional order $\alpha(t)$. \bh{To prevent overfitting, hyper-parameter tuning, specifically for epoch numbers, and an early-stopping mechanism were implemented in our model. The training loss is observed to converge to our desired tolerance while increasing the epoch number.} Our method was implemented using the PyTorch library, and the source code has been made publicly available on GitHub.
\subsection{Fixed Fractional Order FPDE}
In this section, we present the results of experiments using the proposed fPINN-DeepONet framework to solve the fixed fractional-order FPDEs where the fractional order $\alpha_i$ is a real number within $[0,1]$. The fPINN-DeepONet model consists of two sub-networks: the branch net and the trunk net. The branch net encodes the input function $u$ in $100$ fixed sensors, and the trunk net takes the spatiotemporal data [x,t] of size $10,000 \times 1$ as inputs. 
\subsubsection{FODE example}
We consider the simple fixed fractional order FODE example as follow:
\begin{eqnarray}\label{FODE equation}
&  D_t^{\alpha}u(t) = f(t), \qquad  t \in [0,1], \alpha \in [0,1]
\end{eqnarray}
with the exact solution
\begin{eqnarray}\label{FODE equation}
\notag &  u(t) = t^3
\end{eqnarray}
The forcing-term
\begin{eqnarray}
\notag &  f(t)= \frac{\Gamma(4)}{\Gamma(4-\alpha(t))}t^{3-\alpha(t)}
\end{eqnarray}

In this problem, the fractional order $\alpha$ is a real number in [0,1], $\alpha = 0.5$ is applied to solve the above FODE problem. Grid sampling is used to generate time sampling data \bh{$t$} from 100 sensors in the range $[0,1]$. The time sampling data \bh{$t$} serve as the inputs to the branch network, while \bh{$u(t)$} is used as the inputs to the trunk network. We train the model for $100$ epochs, during which the training loss converges to $0.000064$. In addition, we plot the solution to compare the predicted values with the true solution.
\begin{figure}[H]
     \centering
       \includegraphics[width=0.35\textwidth]{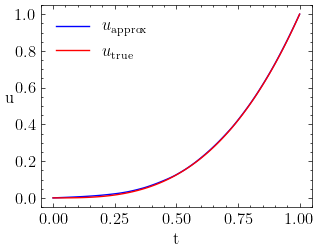}  
     \hfill
     \caption{Fixed Fractional Order FODE Forward Problem: Comparison of the exact solution and fPINN-DeepONet prediction over time.}
    \label{fig:Fixed Fractional Order FODE Forward Problem}
\end{figure}
In Figure~\ref{fig:Fixed Fractional Order FODE Forward Problem},  the plotted curve of the prediction, $u_{approx}$, is observed to closely align with the exact solution, $u_{true}$ for $t \in [0,1]$. According to our calculation, the \bh{$L_2$} relative error is $7.76 \times  10 ^{-5}$, indicating that our model achieves exceptional prediction accuracy.

\subsubsection{FPDE example}
We consider the following fixed fractional-order fractional partial differential equation (FPDE):
\begin{eqnarray}
&&  _0^{C}D_t^{\alpha(x,t)}u(x,t) = \partial_{x}^2 u(x,t)
+f(x,t),\qquad t \in [0,1],~x \in \Omega \equiv [0,\pi],\\
\notag && u(\pi,t) = 0, \quad u(0,t) = 0,\\
\notag && u(x,0) = 0,
\end{eqnarray}
The forcing term \( f(x,t) \) is given by
\begin{eqnarray}
\notag &f(x,t) = \frac{\Gamma(4)}{\Gamma(4 - \alpha(t))} t^{3 - \alpha(t)} \sin(x) + t^3 \sin(x).
\end{eqnarray}
In this problem, the fractional order $\alpha$ is within $[0,1]$ and$\alpha = 0.5$ is employed for this FPDE problem. Spatial sampling data \bh{$x$} is generated using the grid sampling method with $100$ sensors in the range $[0,\pi]$ and time sampling data \bh{$t$} is generated in the range $[0,1]$. The concatenation of spatiotemporal data $[x,t]$ is used as the inputs to the trunk net and $u(x,t)$ is the input function to the branch net. We train the model for $3,000$ epochs as the training loss converges to 0.083. 
\begin{figure}[H]
     \centering
       \includegraphics[width=1.0\textwidth, height=5.0cm]{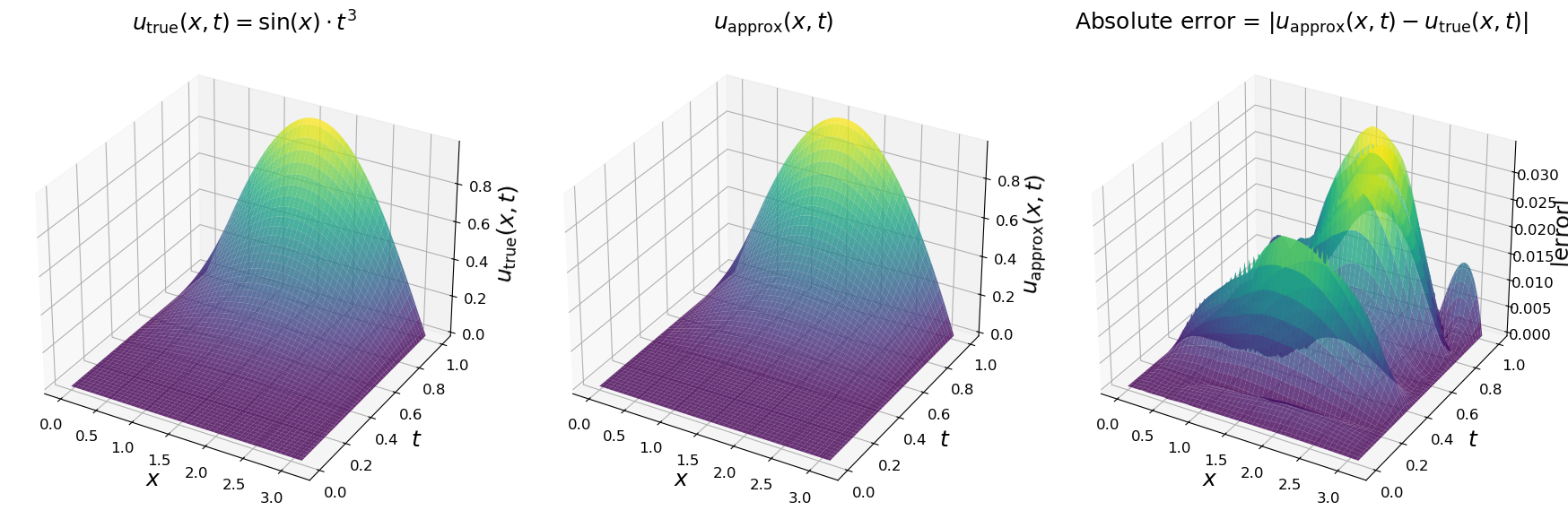}
     \hfill
     \caption{Fixed Fractional Order FPDE Forward Problem: The subfigures (from left to right) display the exact solution, the fPINN-DeepONet prediction, and the corresponding absolute error.}
    \label{fig:Fixed FPDE}
\end{figure}
As shown in Figure~\ref{fig:Fixed FPDE}, the prediction $u_{approx}$ closely matches the true solution $u_{true}$ in both the spatial and temporal domains. The absolute error of 0.0029 and the $L_2$ relative error of $1.75 \times 10^{-5}$ demonstrate the high predictive accuracy of our model.

\subsubsection{\bh{TFMDWEs}}
\bh{In this problem, we study the FPDE with fractional order $\alpha \in [0,1]$ and fractional order $\beta \in [1,2]$.
\begin{eqnarray}\label{One-Dimensional Two-Term FPDE}
&&  _0^{C}D_t^{\alpha(x,t)}u(x,t) = \partial_{x}^2 u(x,t)
+f(x,t),\qquad t \in [0,1],~x \in \Omega \equiv [0,\pi],\\
\notag && u(\pi,t) = 0, \quad u(0,t) = 0,\\
\notag && u(x,0) = 0,
\end{eqnarray}
with the forcing term
\begin{eqnarray}
\notag &    f(x,t) = \frac{\Gamma(4)}{\Gamma(4-\alpha(x,t))}t^{3-\alpha(x,t)} sinx + \frac{\Gamma(4)}{\Gamma(4-\beta(x,t))}t^{3-\beta(x,t)} sinx + t^3 sinx
\end{eqnarray}
} $\alpha = 0.5$ and $\beta = 1.5$ is used to solve this problem. The spatiotemporal sampling data are also generated using grid sampling method with $100$ sensor in the domain. The concatenation of spatiotemporal sampling data $[x,t]$ is used as the input function to the trunk net and $u(x,t)$ is the input function to the branch net. The model is trained with $3,000$ epochs and the training loss converges to $0.0103$. 
\begin{figure}[H]
     \centering
       \includegraphics[width=1.0\textwidth, height=5.0cm]{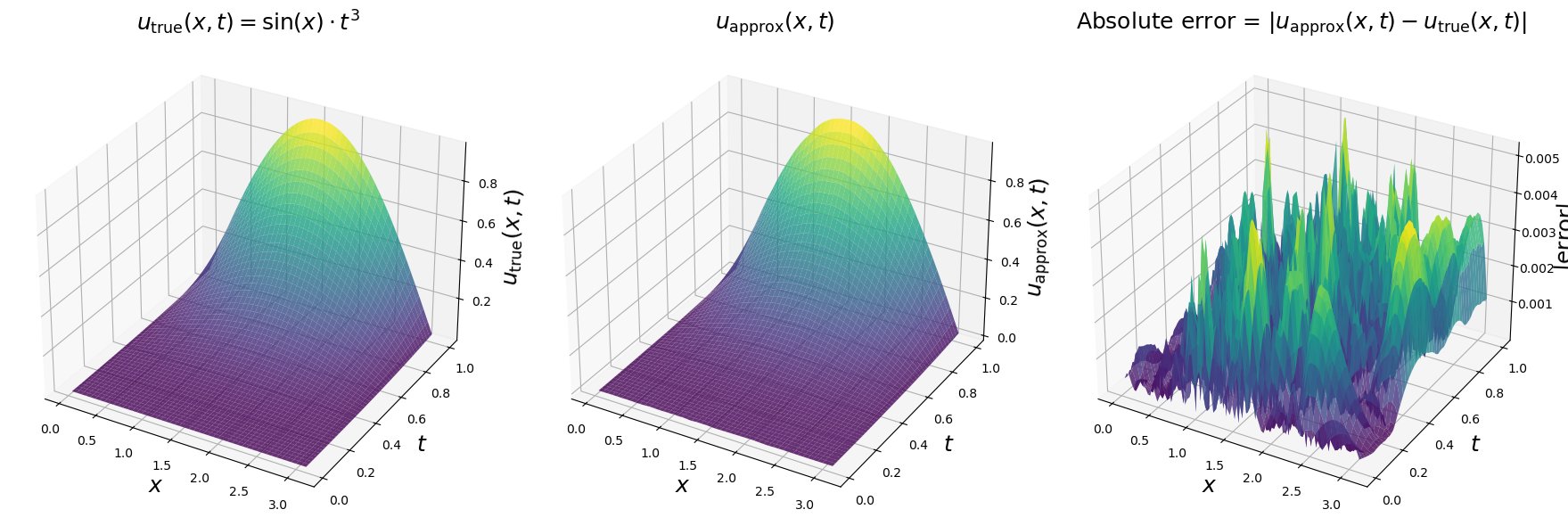}
     \hfill
     \caption{One-Dimensional Two-Term FPDE Forward Problem: The subfigures (from left to right) display the exact solution, the fPINN-DeepONet prediction, and the corresponding absolute error.}
    \label{fig:One-Dimensional Two-Term}
\end{figure}
 Figure~\ref{fig:One-Dimensional Two-Term} illustrates that the predicted solution closely aligns with the true solution in both the spatial and temporal domains. The absolute error is 0.0013 and the $L_2$ relative error is $2.68 \times 10^{-6}$, indicating highly accurate prediction.

\subsection{Variable Fractional-Order FPDE}

In this section, we extend our research to explore the proposed fPINN-DeepONet framework to solve the variable fractional-order FPDEs. Specifically, the fractional order $\alpha_i$ is no longer a real number on $[0,2]$, however, it is a function depending on both time $t$ and space $x$. The fPINN-DeepONet model consists of two subnets: the branch net and the trunk net. The branch net encodes the input function $\alpha_i$ at $100$ fixed sensors, and the trunk net takes the spatio-temporal sampling data [x,t] as the input data.

\subsubsection{FODE Forward example}
The fPINN-DeepONet model is first employed to solve the following FODE problem:
\begin{eqnarray}\label{FODE equation}
&  D_t^{\alpha(t)}u(t) = f(t), \qquad  t \in [0,1]
\end{eqnarray}
with exact solution 
\begin{eqnarray}\label{FODE equation}
\notag u(t) = t^3
\end{eqnarray}
The temporal data $t$ is uniformly sampled in the range $[0,1]$ at $100$ fixed sensors. The fractional order $\alpha(t)$ is a function of $t$ in the range $[0,1].$ The function $u(t)$ was fed as input data to the branch net and the temporal data $t$ was fed as input data to the trunk net.

\begin{figure}[H]
     \centering
     \includegraphics[width=0.35\textwidth]{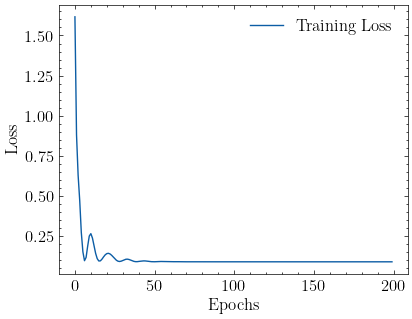}
       \includegraphics[width=0.35\textwidth]{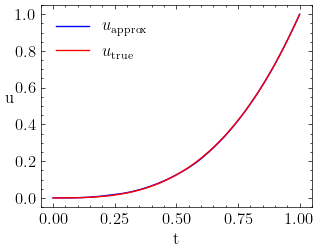}  
     \hfill
     \caption{Variable Fractional Order FODE Forward Problem: The subfigure on the left shows the training loss over epochs, while the subfigure on the right presents a comparison between the exact solution and the fPINN-DeepONet prediction over time.}
    \label{fig:Variable Fractional Order FODE Forward Problem}
\end{figure}
The FODE results using the fPINN-DeepONet model are presented in Figure~\ref{fig:Variable Fractional Order FODE Forward Problem}, It is evident that an excellent agreement can be achieved between the prediction of fPINN-DeepONet and the ground truth. During the network training, the training loss converges  to our desired tolerance after training for $200$ epochs. The subfigure on the right shows the comparison between the predicted and exact solutions for a random test input sample. To better visualize the difference, the relative $L_2$ error is calculated to be $1.3892 \times 10^{-5}$, which demonstrates the prediction accuracy of our proposed model.
 
\subsubsection{FODE Inverse example}

In this section, the capability of the fPINN-DeepONet model is examined in solving inverse problems. Given the true solution \( u(t) \) and the right-hand-side function \( f(t) \), our objective is to approximate the fractional order \( \alpha(t) \). Specifically, \( \alpha(t) \) is used as the input to the branch network, while the temporal variable \( t \), defined within the domain of \( \alpha \), serves as the input to the trunk network.

\begin{figure}[H]
     \centering
      \includegraphics[width=0.35\textwidth]{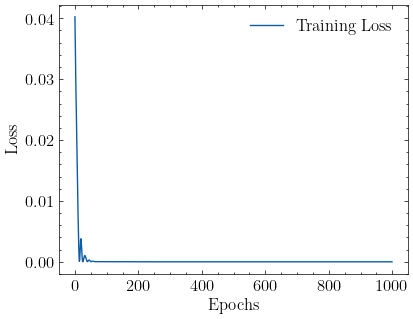}  
       \includegraphics[width=0.35\textwidth]{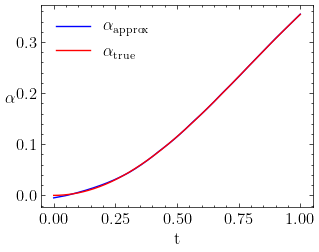}  
     \hfill
     \caption{Variable Fractional Order FODE Inverse Problem: The subfigure on the left shows the training loss over epochs, while the subfigure on the right compares the exact fractional order with the fPINN-DeepONet prediction over time.}
    \label{fig:Variable Fractional Order FODE inverse Problem}
\end{figure}
As shown in Figure~\ref{fig:Variable Fractional Order FODE inverse Problem}, the subfigure on the left shows the training loss using the fPINN-DeepONet model. The loss value converges to the desired tolerance for approximately 500 iterations. The subfigure on the right shows the comparison between the approximated fractional order value and the real value. To better understand the difference, we calculated the \bh{$L_2$} relative error $9.2195 \times 10^{-5}$ showing the excellent prediction of our model on the inverse problem.
\subsection{FPDE example}
\label{FPDE example}
In this section, the fPINN-DeepONet model is examined on the FPDE problem:
\begin{eqnarray}\label{multi-dw}
&&  _0^{C}D_t^{\alpha(t)}u(x,t) = \partial_{x}^2 u(x,t)
+f(x,t),\qquad t \in [0,1],~x \in \Omega \equiv [0,\pi],\\
\notag && u(\pi,t) = 0, \quad u(0,t) = 0,\\
\notag && u(x,0) = 0,
 \label{FPDE equation}
\end{eqnarray}
The forcing term \( f(x,t) \) is given by
\begin{eqnarray}
\notag &f(x,t) = \frac{\Gamma(4)}{\Gamma(4 - \alpha(t))} t^{3 - \alpha(t)} \sin(x) + t^3 \sin(x).
\end{eqnarray}
The time data $t$ and the spatial data $x$ were both sampled using the grid sampling method. The time data $t$ are located within [0,1] and the spatial data $x$ are located on [0,$\pi$]. The fractional order $\alpha(t)$ is a function of $t$ in the field [0,1].
\subsubsection{FPDE Forward example}
Given the right-had side ~$f(x,t) = \frac{\scriptstyle\Gamma(4)}{\scriptstyle\Gamma(4-\alpha(t))}\cdot t^{{3-\alpha(t)}} \cdot sin(x)+t^3 \cdot sin(x) $ and the fractional order $\alpha(t) = \frac{1}{2}sin^2 t$, the fPINN-DeepONet was constructed to approximate the solution $u(x,t)$ of the corresponding FPDE. In this model, the input to the branch net is $\alpha(t)$, while the trunk net takes $[x,t]$ as input. The residual loss is constructed as:
\begin{eqnarray}\label{FPDE loss}
Loss =  \int_{0}^1  \int_{0}^{\pi} |D_t^{\alpha_{NN}(t)}u(x,t)- \frac{\partial^2 u}{\partial x^2} - f(t) | dxdt
\end{eqnarray}
    \begin{figure}[H]
     \centering
        
       \includegraphics[width=1.0\textwidth, height=5.0cm]{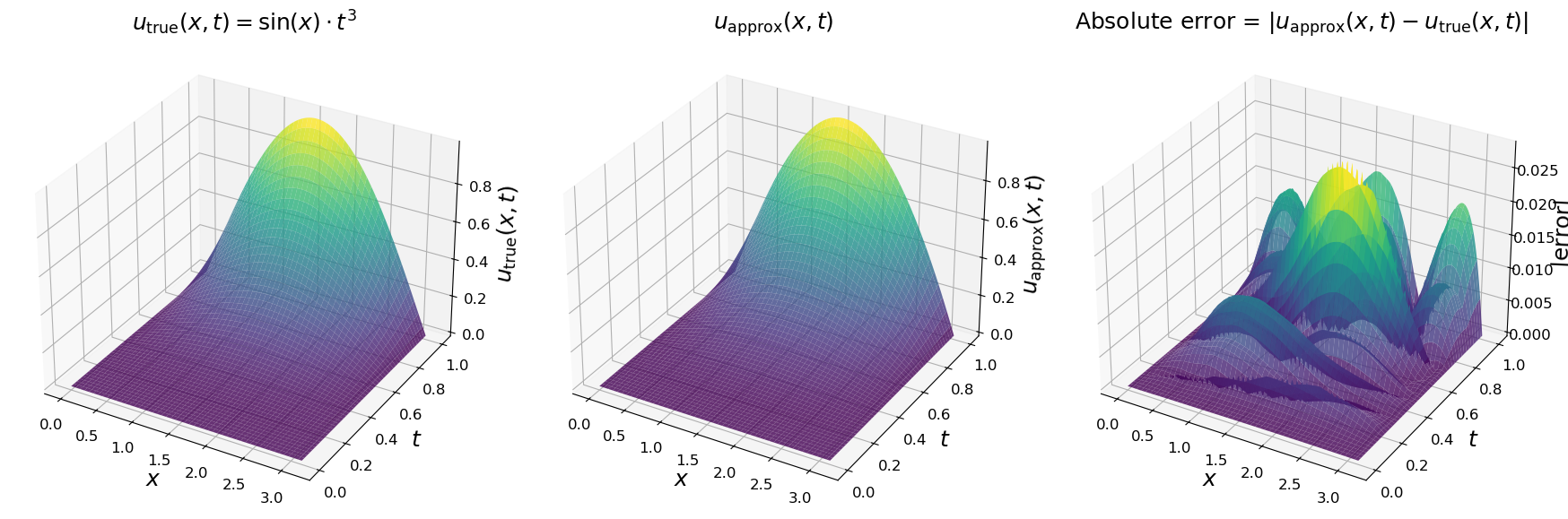}

     \hfill
     \caption{FPDE Forward Problem: The subfigures (from left to right) display the exact solution, the fPINN-DeepONet prediction, and the corresponding absolute error.}
    \label{fig:FPDE Forward Problem}
\end{figure}
As shown in Figure~\ref{fig:FPDE Forward Problem}, the prediction value $u_{approx}(x,t)$ closely matches the ground truth $u_{truth}$ in the spatial-temporal domain. The $L_2$ relative error is calculated to be 0.0001, which is consistent with the visual agreement observed in the figure.

\subsubsection{FPDE Inverse example}
To evaluate the capability of fPINN-DeepONet in solving inverse problems for FPDEs, we constructed a model to approximate the unknown fractional order $\alpha(t)$. Given the right-hand side $f(x,t)$ and information of the real solution $u(x,t) = t^3sinx$, the fPINN-DeepONet framework was employed, with the branch network designed to learn $\alpha(t)$ and the trunk network to represent the solution u(x,t). The loss function used in this inverse problem setup is identical to that defined for the FPDE forward problem~\ref{FPDE loss}.
\begin{figure}[H]
     \centering
     \includegraphics[width=0.35\textwidth]{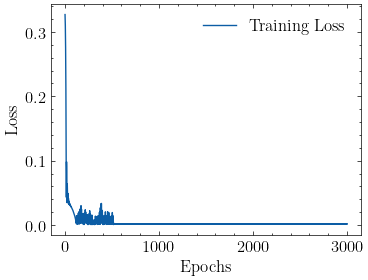}
       \includegraphics[width=0.35\textwidth]{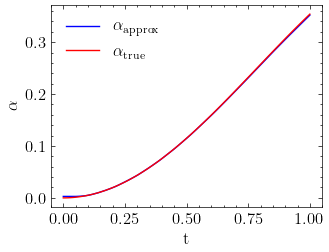}

     \hfill
     \caption{FPDE Inverse Problem: The subfigure on the left shows the training loss as a function of epochs, while the subfigure on the right presents a comparison between the true fractional order and the fPINN-DeepONet predicted fractional order over time.}
    \label{fig:FPDE Inverse Problem}
\end{figure}
As shown in Figure~\ref{fig:FPDE Inverse Problem}, the training loss converges to the desired tolerance within approximately 500 iterations. The figure on the right compares the predicted and true values of the fractional order, illustrating that the approximation closely matches the true \(\alpha(t)\) across the domain \([0,1]\). To further quantify the model’s accuracy, the \(L_2\) relative error was computed to be \(8.8114 \times 10^{-6}\), confirming the high precision of the proposed approach.

\subsection{FPDE With Fractional Order $\alpha(x,t)$}
In this section, we extend the fractional order \(\alpha(t)\) to a function of both time \(t\) and space \(x\), denoted as \(\alpha(x,t)\).
 We consider the following equation: 
\begin{eqnarray}\label{multi-dw}
&&  _0^{C}D_t^{\alpha(x,t)}u(x,t) = \partial_{x}^2 u(x,t)
+f(x,t),\qquad t \in [0,1],~x \in \Omega \equiv [0,\pi],\\
\notag && u(\pi,t) = 0, \quad u(0,t) = 0,\\
\notag && u(x,0) = 0,
 \label{FPDE equation with alpha(x,t)}
\end{eqnarray}
Where the fractional order $\alpha(x,t)$ is within [0,1].
\subsubsection{Forward example}
Given the right-hand side  
\[
f(x,t) = \frac{\Gamma(4)}{\Gamma(4 - \alpha(x,t))} \cdot t^{3 - \alpha(x,t)} \cdot \sin(x) + t^3 \cdot \sin(x)
\]  
and the fractional order  
\[
\alpha(x,t) = \frac{1}{2} \sin(x) \cdot t^2,
\]  
we construct the fPINN-DeepONet model to approximate the true solution \( u(x,t) \) of the FPDE. In this setup, the branch network processes the input term \(\alpha(x,t)\), while the trunk network handles the input \([x, t]\) to learn the solution \(u(x,t)\). The loss function used is the same as that defined for the FPDE forward problem (see Equation~\ref{FPDE loss}).

 \begin{figure}[H]
     \centering
       \includegraphics[width=1.0\textwidth, height=5.0cm]{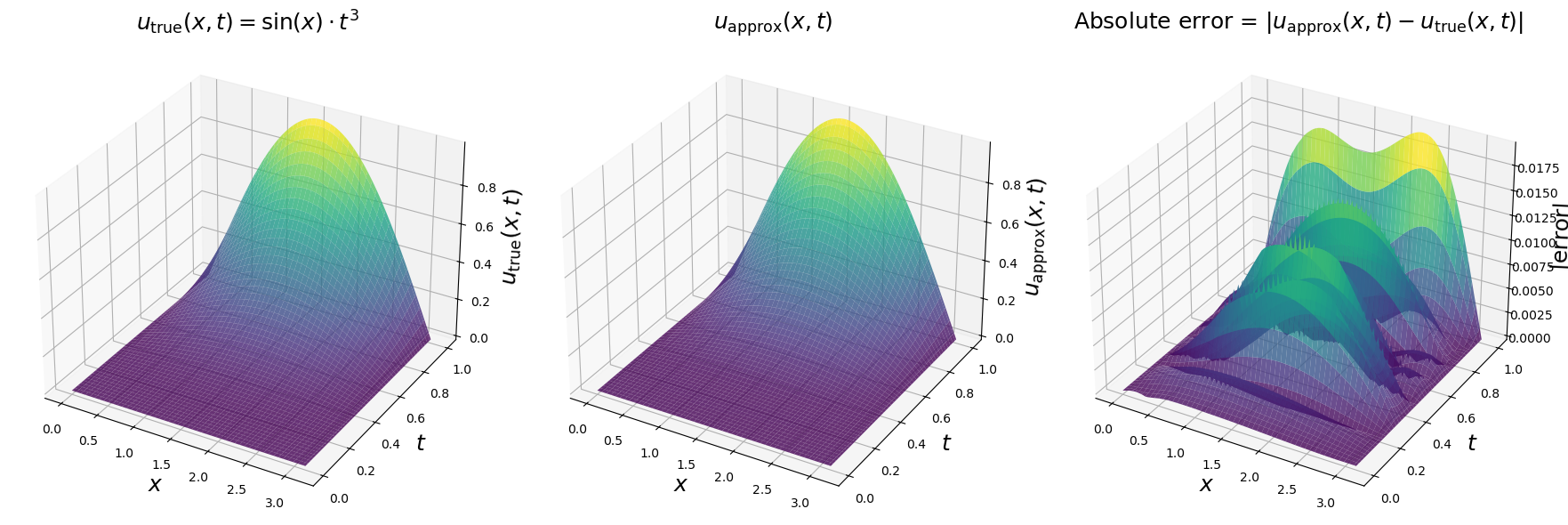}
     \hfill
     \caption{$\alpha(x,t)$ FPDE Forward Problem: The subfigures (from left to right) display the exact solution, the fPINN-DeepONet prediction, and the corresponding absolute error.}
    \label{fig:FPDE alpha(x,t) Forward Problem}
\end{figure}
As shown in Figure~\ref{fig:FPDE alpha(x,t) Forward Problem}, the predicted solution closely matches the real solution $u(x,t)$. The \(L_2\) relative error was computed to be \(5.5107 \times 10^{-5}\), demonstrating the high precision of the model.

\subsubsection{Inverse example}
This section investigates the fractional-order function \(\alpha(x,t)\) using the fPINN-DeepONet framework. Given the exact solution \(u(x,t) = t^3 \sin(x)\) and the corresponding right-hand side \(f(x,t)\), the model is constructed to approximate the underlying fractional order \(\alpha(x,t)\). The true fractional-order function is \(\alpha(x,t) = \frac{1}{2} \sin(x) t^2\). The \(L_2\) relative error is calculated as \(3.1652 \times 10^{-5}\), indicating the high accuracy of the approximation.

\begin{figure}[H]
     \centering
       \includegraphics[width=1.0\textwidth, height=5.0cm]{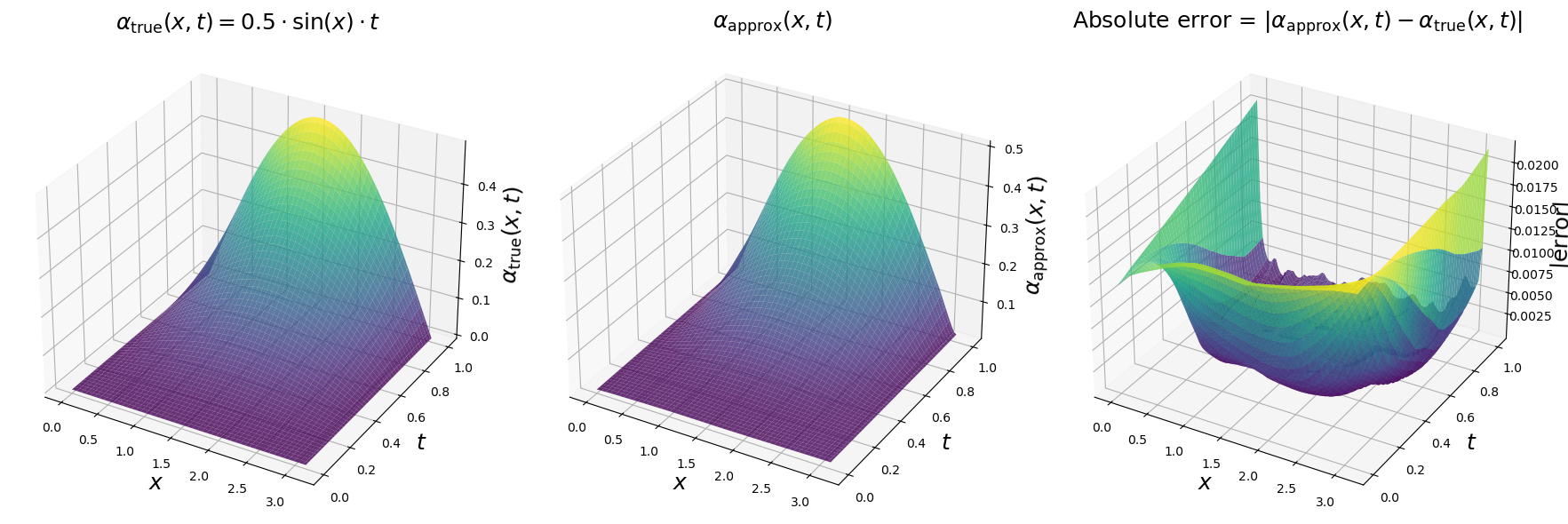}
     \hfill
     \caption{$\alpha(x,t)$ FPDE Inverse Problem: The subfigures (from left to right) show the exact fractional order, the fPINN-DeepONet prediction, and the corresponding absolute error.}
    \label{fig:FPDE alpha(x,t) Inverse Problem}
\end{figure}
\subsection{\bh{High-dimensional FPDE}}
\bh{
    In this section, we extended the spatial domain to multiple dimensions, e.g., $(x_1, x_2,\dots x_n,t)$ instead of just $x,t$ by introducing another spatial variable $x_i \in [0,\pi]$. The equation can be shown as:
\begin{eqnarray}\label{multi-dw}
_0^{C}D_t^{\alpha(t)}u(x_1, x_2, \dots, x_n, t) &=& \sum_{i=1}^{d} \frac{\partial^2 u}{\partial x_i^2} 
+ f(x_1, x_2, \dots, x_n, t), \nonumber \\ 
\end{eqnarray}
Where \begin{eqnarray}
    && \quad t \in [0,1], \quad \mathbf{x} = (x_1, x_2, ..., x_d) \in \Omega \equiv [0, \pi]^d.\nonumber
\end{eqnarray}
The forcing term: 
 \begin{eqnarray}
\notag f(x_1, x_2, \dots, x_d, t) & = & \left( \frac{\Gamma(4)}{\Gamma(4 - \alpha)} \right) t^{3 - \alpha} \prod_{i=1}^{d} \sin(x_i) + t^3 \prod_{i=1}^{d} \sin(x_i)
\end{eqnarray}
In this experiment, an additional spatial variable \(y\) is introduced. Both spatial variables \(x\) and \(y\) are uniformly sampled from the interval \([0,\pi]\), with 100 data points in each dimension. Similarly, the temporal variable \(t\) is sampled from \([0,1]\) using 100 data points. The trunk network receives \([x, y, t]\) as input, while the branch network takes \(\alpha(t)\) as input. A total of 10{,}000 data points are allocated for the initial and boundary conditions, and 5{,}000 data points are used for FPDE collocation points. The model is trained for 1{,}000 epochs until the training loss fully converges. The final \(L_2\) loss between the predicted and true solutions is 0.0499.
\begin{figure}[H]
     \centering
       \includegraphics[width=0.45\textwidth]{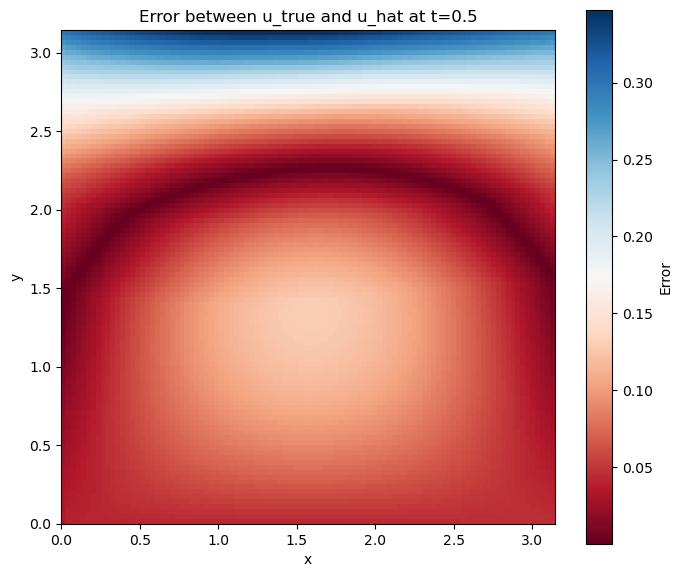} 
     \hfill
     \caption{Heatmap of testing error of high-dimensional FPDE at $t=0.5.$}
    \label{fig:Heatmap of the Error of High-dimensional FPDE}
\end{figure}
Figure~\ref{fig:Heatmap of the Error of High-dimensional FPDE} shows the heatmap of testing error when $t=0.5$ lower part of figure are colored as red which correspond to error of 0.05 and top part of figure are white to blue which correspond to error of 0.20, which is consistent with our previous calculation and proves the feasibility of using  our proposed  to solve high-dimensional FPDEs. 
}

\subsection{\bh{Robustness}}
 \bh{This section evaluates the robustness of the proposed model by introducing noise into the training data. The FPDE forward problem is defined in Eq.35, where the temporal data t is within [0,1] and the spatial data x is within [0,$\pi$]. The fractional order $\alpha(t)$ is treated as a function of t over the same interval. The branch network receives $\alpha(t)$ as input, while the trunk network processes the solution $u(x,t)$. To simulate noisy conditions, a noise term is added to the right-hand side forcing term $f(x,t)$:
   \begin{eqnarray}\label{Noise equation}
&  f(x,t) = \hat{f}(x,t) + noise
\end{eqnarray}
Gaussian noise of varying intensity is introduced into the forcing term, and the resulting test errors under different noise levels are summarized in Table~\ref{tab:noise_comparison}.
\begin{table}[H]
    \centering
    \begin{tabular}{lccc}
        \toprule
        Noise Level & 20\% data noise & 50\% data noise & 80\% data noise \\
        \midrule
        Testing Error & 0.0245 & 0.0268 & 0.0267 \\
        \bottomrule
    \end{tabular}
    \caption{Testing error of fPINN-DeepONet under different noise levels.}
    \label{tab:noise_comparison}
\end{table}
As we can observe from Table~\ref{tab:noise_comparison}, the testing error of our proposed model falls in  $\sim 10^{-2}$ with Gaussian noise range from $20\%$ to $80\%$. 
Also, we plotted the true, predicted solutions and the absolute error shown the figure of error:
\begin{figure}[htbp]
    \centering
    \begin{minipage}{0.32\textwidth}
        \centering
        \includegraphics[width=\linewidth]{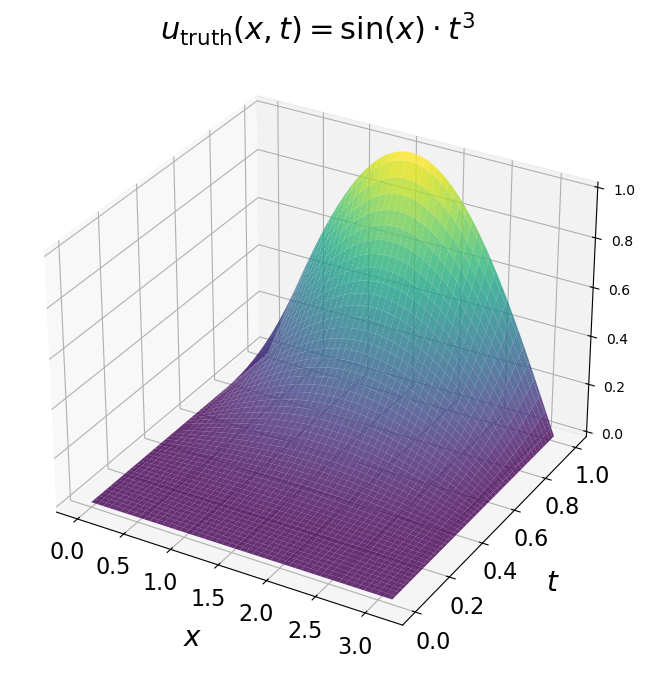}
    \end{minipage}
    \hfill
    \begin{minipage}{0.32\textwidth}
        \centering
        \includegraphics[width=\linewidth]{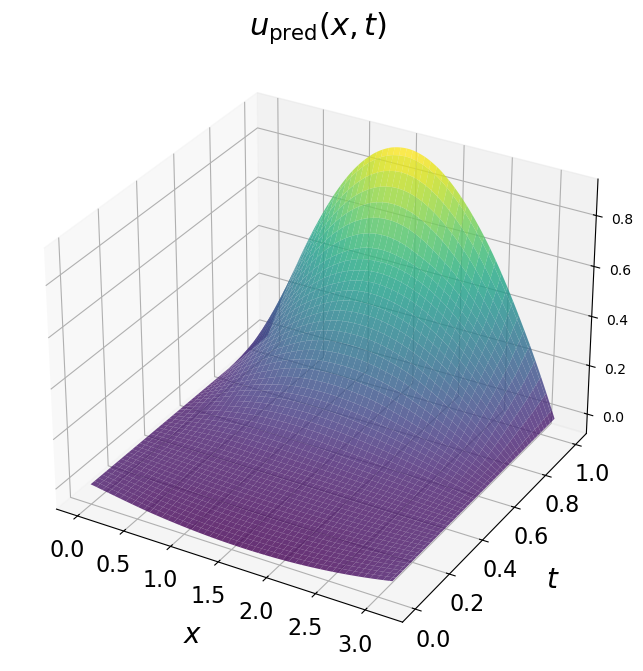}
    \end{minipage}
    \hfill
    \begin{minipage}{0.32\textwidth}
        \centering
        \includegraphics[width=\linewidth]{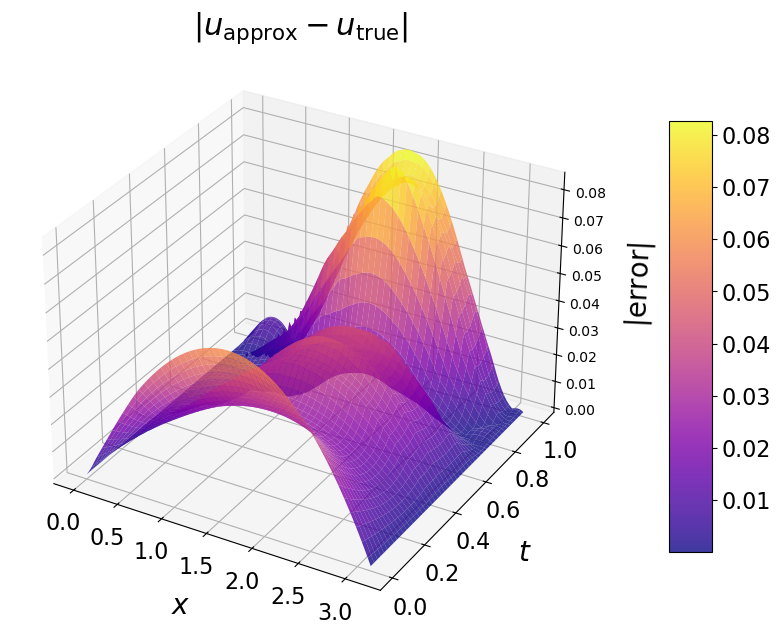}
    \end{minipage}
    \caption{Subfigures (from left to right) display the exact solution, the fPINN-DeepONet prediction, and the absolute error under $50\%$ Gaussian noise.}
    \label{fig:noisy data}
\end{figure}
The figures we shown in figure~\ref{fig:noisy data} show the true, predicted solutions and the absolute error between the true solution and the predicted solution under $50\%$ Gaussian noise. From our observation, the result is consistent with the testing error table that the final testing loss falls in $\approx 10^{-2}$, which demonstrated the robustness of our proposed model with respect to data noise.}

      
       

\subsubsection{Discontinuous fractional order}

This section investigates the solution of an FODE with a discontinuous fractional order. The fractional order is defined in Equation~\ref{Discontinuous fractional order equation} as follows:

\begin{equation}\label{Discontinuous fractional order equation}
\alpha_{\text{true}}(t) = 
\begin{cases} 
    t & \text{for } 0 \leq t < 0.5 \\
    0.9 & \text{for } 0.5 \leq t < 1 \\
    0 & \text{otherwise}
\end{cases}
\end{equation}

The fractional order function exhibits discontinuities at \(t = 0\), \(t = 0.5\), and \(t = 1\). Temporal data \(t\) is sampled from the interval \([0,1]\), with 100 data points. The fPINN-DeepONet framework is constructed with \(\alpha(t)\) as the input to the branch network and \(t\) to the trunk network. The model is trained for $8,000$ epochs until convergence, with the training loss reaching 0.053.

\begin{figure}[H]
     \centering

        \includegraphics[width=0.35\textwidth]{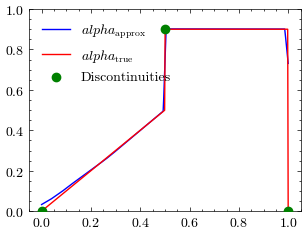}
     \hfill
     \caption{Discontinuous fractional order Inverse Problem: A comparison between the true fractional order and the predicted fractional order obtained from the fPINN-DeepONet model.}
    \label{fig:Discontinuous fractional order result}
\end{figure}

Figure~\ref{fig:Discontinuous fractional order result} shows the comparison between the true and predicted fractional order values. The discontinuities at \(t = 0\), \(t = 0.5\), and \(t = 1\) are clearly observed. The predicted fractional order, \(\alpha_{\text{approx}}\), closely matches the true value \(\alpha_{\text{true}}\). The \(L_2\) relative error is calculated as 0.0014, confirming the accuracy of the fPINN-DeepONet model.

\section{Hyperparameter Testing}
To obtain the best performance of the fPINN-DeepONet model, we investigated the variation in test error with increasing three hyperparameters: the number of data points, the width of the network, and the depth of the network in the FPDE example.
\begin{figure}[H]
\centering
\includegraphics[width=0.3\textwidth]{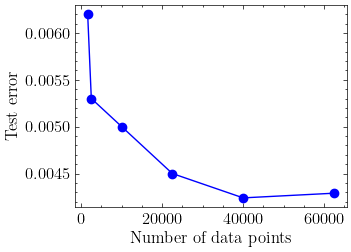}
\includegraphics[width=0.3\textwidth]{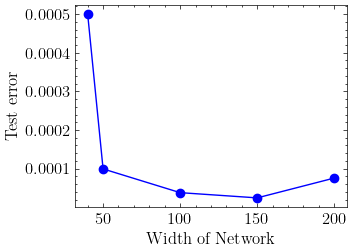}
\includegraphics[width=0.3\textwidth]{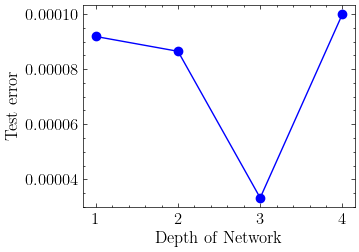}
\caption{Hyperparameter Testing: The subfigure on the left illustrates test error varies with an increasing number of data points. The center subfigure shows the variation in test error with increasing network width, while the subfigure on the right depicts the change in test error as network depth increases.}
\label{fig:hyperparameter_test}
\end{figure}

As shown in Figure~\ref{fig:hyperparameter_test}, the subfigure on the left presents the test results for varying numbers of data points. We set the neural network width to 100 and depth to 2. The results indicate that when using $200 \times 200$ data points for training, the test error reaches approximately $10^{-5}$. The right two subfigures in Figure~\ref{fig:hyperparameter_test} illustrate the effect of varying the width and depth of the neural network. Among all tested hyperparameters, the smallest test error, around $10^{-5}$, is achieved with a network width of 150, depth of 3, and $200 \times 200$ data points.

\section{Discussion and Conclusion}
The primary objective of our work is to apply physics-informed deep operator learning techniques to solve fractional partial differential equations (FPDEs). We address several previously unexplored challenges in the application of deep learning to FPDEs. Our first contribution is the derivation of an $L_2$ approximation for solving FPDEs with fractional orders $\alpha \in [0, 2]$, which we incorporate into the proposed fPINN-DeepONet framework to formulate the residual loss function. Second, we investigate two major classes of FPDEs: fixed fractional order and variable fractional order. To address these, we develop two distinct algorithms within the fPINN-DeepONet framework. For fixed fractional orders, the input function $u(x,t)$ is mapped to the output; for variable fractional orders, the input function $\alpha(x,t)$ is used to generate the output. Third, we explore a variety of FPDE problems, including fractional-order ordinary differential equations (FODEs), forward and inverse FPDE problems, spatial-temporal fractional orders, and discontinuous fractional orders. Finally, we assess the model’s performance under noisy data conditions, demonstrating that the proposed framework not only achieves high prediction accuracy but also reduces data requirements and remains robust in challenging noisy environments.

\section*{Declaration of competing interest}
The authors declare that they have no known competing financial interests or personal relationships that could have appeared to
influence the work reported in this paper.

\section*{Data availability}
\url{https://github.com/jeremylu916/fPINN-DeepONet}

\section*{Acknowledgment}
The authors appreciate Prof. Zhi-zhong Sun for providing us useful suggestions through careful reading.

G. Lin would like to thank the support by National Science Foundation (DMS-2053746, DMS-2134209, ECCS-2328241, CBET-2347401 and OAC-2311848), and U.S.~Department of Energy (DOE) Office of Science Advanced Scientific Computing Research program DE-SC0023161, and DOE–Fusion Energy Science, under grant number: DE-SC0024583. 

\bibliographystyle{abbrv}


\end{document}